\documentclass[a4paper]{amsart}

\usepackage{amsthm}
\usepackage{amsfonts}
\usepackage{amsmath,amsthm,amssymb}
\usepackage{comment}
\usepackage{fancybox}
\usepackage{epsf}
\usepackage[dvips]{graphicx}

\newtheorem{theorem}{Theorem}[section]
\newtheorem{lemma}[theorem]{Lemma}

\theoremstyle{definition}

\theoremstyle{remark}
\newtheorem{remark}[theorem]{Remark}
\newtheorem{example}[theorem]{Example}
\newtheorem{claim}[theorem]{Claim} 
\newtheorem{observation}[theorem]{Observation}

\numberwithin{equation}{section}
\begin{document} 

\title[Milnor invariants of length $2k+2$ for links]{Milnor invariants of length $2k+2$ 
for links with vanishing Milnor invariants of length $\leq k$}

\author[Yuka Kotorii and Akira Yasuhara]{Yuka Kotorii \\ Akira Yasuhara}

\begin{abstract} 
J.-B. Meilhan and the second author showed that any Milnor $\overline{\mu}$-invariant of length 
between $3$ and $2k+1$ can be represented as a combination of  HOMFLYPT polynomial of knots obtained by 
certain band sum of the link components, if all $\overline{\mu}$-invariants 
of length $\leq k$ vanish.  
They also showed that their formula does not hold for length $2k+2$.  
In this paper, we improve their formula to give the $\overline{\mu}$-invariants of length $2k+2$ 
by adding correction terms. The correction terms can be given by a 
combination of  HOMFLYPT polynomial of knots determined by $\overline{\mu}$-invariants  of length $k+1$.  
In particular, for any 4-component link the $\overline{\mu}$-invariants of length $4$ are given by our formula, since 
all $\overline{\mu}$-invariants of length $1$ vanish. 

\end{abstract} 

\thanks{ The second author is partially supported by a Grant-in-Aid for Scientific Research (C) 
($\#$23540074) of the Japan Society for the Promotion of Science.}

% \date{\today}

\maketitle

%%%%%%%%%%%%%%%%%%%%%%%%%%%%%%%%%%%%%%%%%%%%%%%%%%%%%%%%%%%%%%%%%%%%%%%%%%%%%%%%%%%%%%%%%%%%%%%%%%%%%%%%%%%%%%%%%%%%%%%%%%%%%%%%%%
\section{Introduction}
For an ordered, oriented link in the 3-sphere, J. Milnor \cite{Milnor, Milnor2} defined a family of invariants, 
known as {\em Milnor $\overline{\mu}$-invariants}. 
For an $n$-component link $L$, Milnor invariant is specified by a sequence $I$ of numbers in $\{1,2,\ldots ,n\}$
and denoted by $\overline{\mu}_L(I)$. The length of the sequence $I$ is called the {\em length} 
of the Milnor invariant $\overline{\mu}_L(I)$. 
It is known that Milnor invariants of length two are just linking numbers.
In general, Milnor invariant $\overline{\mu}_L(I)$ is only well-defined modulo the greatest common divisor $\Delta _L(I)$ 
of all Milnor invariants $\overline{\mu}_L(J)$ such that $J$ is obtained from $I$ by removing at least one 
index and permuting the remaining indices cyclicly. 
%This indeterminacy comes from the choice of the meridian curves generating the link group.
If the sequence is non-repeated, then this invariant is also link-homotopy invariant and 
we call it {\em Milnor link-homotopy invariant}.
Here, the {\em link-homotopy} is an equivalence relation generated by self-crossing changes.   

In \cite{P}, M. Polyak gave a formula expressing Milnor invariant of length 3, 
and in \cite{MY}, J-B. Meilhan and the second author generalized it. 
More precisely, in \cite{MY} they showed that any Milnor invariant of length between $3$ and 
$2k+1$ can be represented as a combination of  HOMFLYPT polynomial of knots obtained by certain band 
sum of the link components, if all Milnor invariants of length $\leq k$ vanish.  
Their assumption that a link has vanishing Milnor invariants of length $\leq k$ is essential 
to compute Milnor invariants of length up to $2k+1$ via their formula. In fact, their formula 
does not hold for length $2k+2$ (\cite[Section 7]{MY}). 

In  this paper, we 
improve their formula to give the Milnor invariants of length $2k+2$ 
by adding correction terms. 
Our formula implies that any Milnor invariant of length $2k+2$ can be given by 
a combination of  HOMFLYPT polynomial of knots obtained by certain band 
sum operations and knots determined by the first non vanishing Milnor invariants, 
which are Milnor invariants of length $k+1$ (Theorem~\ref{main}). 
In particular, the Milnor invariants of length $4$ for any link are given by our formula, since 
all Milnor invariants of length $1$ vanish by the definition (Theorem~\ref{main2}).

Recall that the HOMFLYPT polynomial of a knot $K$ is of the form $P(K;t,z)=\sum_{k=0}^{N} P_{2k}(K;t)z^{2k}$ 
and denote by $P^{(l)}_0(K)$ the $l$-th derivative of $P_0(K;t) \in \mathbb{Z}[t^{\pm}]$ evaluated at $t=1$.
Denote by $(\log P_0(K))^{(l)}$ the $l$-th derivative of $\log P_0(K;t)$ evaluated at $t=1$.
We note that $\log P_0(K;t)$  is an  additive invariant 
for knots under the connected sum, 
since  the HOMFLYPT polynomial of knots is multiplicative. 
In particular, $(\log P_0(K))^{(l)}$ is additive.  
It is known that $P_{0}^{(l)}$ is a finite type invariant of degree $l$ \cite{KM}.  
Since $(\log P_0(K))^{(l)}$ is equal to $P_0(K)^{(l)}$ plus a sum of products of $P_0(K)^{(k)}$'s with $k<l$,  
$(\log P_0)^{(l)}$ is an additive finite type knot invariant of degree $l$. 
We also notice that 
$(\log P_0(K))^{(l)} = P^{(l)}_0(K)$ for $l=1,2,3$, since $P_0(K;1)=1$ and $P_0^{(1)}(K;1)=0$.

Let $L=\bigcup_{i=1}^n L_i$ be an $n$-component link in $S^3$.
Let $I=i_1i_2\ldots i_{m}$ be a sequence of $m$ distinct elements of $\{1,2,\ldots ,n\}$.  
Let $B_I$ be an oriented $2m$-gon, and let $p_j~(j=1,2,\ldots ,m)$ denote $m$ mutually disjoint
 edges of $\partial B_I$ according to the boundary orientation. 
Suppose that $B_I$ is embedded in $S^3$ such that $B_I\cap L=\bigcup_{j=1}^{m} p_j$, and 
such that each $p_j$ is contained in $L_{i_j}$ %the ${i_j}$-th component of $L$ 
with opposite orientation. 
We call such a disk an \emph{I-fusion disk} for $L$.  
For any subsequence $J$ of $I$, 
we define the oriented knot $L_J$ as the closure of  
$\left( (\bigcup_{i\in \{J\}} L_i)\cup \partial B_I \right)\setminus \left( (\bigcup_{i\in \{J\}} L_i)\cap \partial B_I \right)$,  
where $\{J\}$ is the subset of $\{1,2,\ldots ,n\}$ formed by all indices appearing in the sequence $J$. 

Given a sequence $I$ of elements of $\{1,2,\ldots ,n\}$, the notation 
$J<I$ will be used for any subsequence $J$ of $I$, possibly empty or equal to $I$ itself, and  
$|I|$ will denote the length of the sequence $I$.

\medskip
\begin{theorem}\label{main}
Let $L$ be an $n$-component link in $S^3$ ($n \geq 4$) with vanishing Milnor link-homotopy 
invariants of length $\leq k$. 
Then for any sequence $I$ of length $2k+2$ of elements of $\{1,2, \ldots ,n \}$ without repeated
 number and for any $I$-fusion disk for $L$, we have 
\[ \overline{\mu}_{L}(I) \equiv - \frac{1}{(2k+1)!2^{2k+1}} \sum_{J<I} (-1)^{|J|}(\log P_0(L_J))^{(2k+1)}
 - \delta_L(I) \pmod {\Delta _L(I)},  \]
where 
$\delta_L(I)$ is an invariant of $L$  that determined by Milnor invariants for length-$(k+1)$ 
subsequences of $I$ which is defined in Subsection 2.5.
\end{theorem}

\noindent
With the same assumption as in Theorem \ref{main}, the same formula but $\delta_L(I)=0$ 
holds for a sequence $I$ with $3 \leq |I| \leq 2k+1$ \cite{MY}.  

We also give the case of 4-component links more clearly. 

\medskip
\begin{theorem}\label{main2}
Let $L$ be a 4-component link in $S^3$.
Then for any sequence $I=i_1i_2i_3i_4$ of distinct elements of $\{1,2,3,4\}$ and for any $I$-fusion disk for $L$, we have 
\[ \overline{\mu}_{L}(I) \equiv - \frac{1}{48} \sum_{J<I} (-1)^{|J|}P^{(3)}_{0}(L_J)
-\frac{1}{2} x_{i_1i_3}x_{i_2i_4}(x_{i_1i_3} +x_{i_2i_4} -1)  \pmod {\Delta _L(I)},  \]
where $x_{ij}$ is the linking number of $i$-th component and  $j$-th component of $L$.
\end{theorem}

\medskip
\begin{remark}
We note that $x_{ij}$ is divisible by $\Delta_L(I)$ if $ij$ is a subsequence of $I$. 
Hence the correction term 
$\frac{1}{2} x_{i_1i_3}x_{i_2i_4}(x_{i_1i_3} +x_{i_2i_4} -1)$  
vanishes up to modulo $\Delta_L(I)$ if either $x_{i_1i_3}$ or $x_{i_2i_4}$ is even.
\end{remark}

\medskip
\begin{remark}
We can generalize Theorem \ref{main} and Theorem \ref{main2} about all repeated sequences 
by the same arguments as those in \cite[Introduction]{MY}. 
That is, we have formulae for not only Milnor link-homotopy invariants but also Milnor isotopy invariants.   
\end{remark}

%%%%%%%%%%%%%%%%%%%%%%%%%%%%%%%%%%%%%%%%%%%%%%%%%%%%%%%%%%%%%%%%%%%%%%
%%%%%%%%%%%%%%%%%%%%%%%%%%%%%%%%%%%%%%%%%%%%%%%%%%%%%%%%%%%%%%%%%%%%%%

\bigskip
\section{Preliminary}

\medskip
\subsection{String link} 
Let $n$ be a positive integer, and let $D^2 \subset  \mathbb{R}^2$ be the unit disk equipped with 
$n$ marked points $x_1, x_2, \ldots, x_n$ in its interior, lying in the diameter on the $x$-axis of  $\mathbb{R}^2$. 
An {\em $n$-string link} (or $n$-component string link) is the image of a proper embedding 
$\sqcup_{i=1}^n[0,1]_i \rightarrow D^2 \times [0,1]$ of the disjoint union $\sqcup_{i=1}^n[0,1]_i $ 
of $n$ copies of $[0,1]$ in $D^2 \times [0,1]$, such that for
each $i$ the image of $[0,1]_i $ runs from $(x_i,0)$ to $(x_i,1)$. 
Each string of an $n$-string link is equipped with an (upward) orientation. 
The $n$-string link $\{x_1,x_2,\ldots,x_n\}\times [0,1]$ in $D^2 \times [0,1]$ is called the {\em trivial $n$-string link} 
and denoted by ${\bf 1}_n$. 
Let $y_1,y_2,\ldots ,y_n$ be points in $\partial D^2$ in Figure~\ref{fig:disk},
$p_i=x_iy_i~(i=1,2,\ldots ,n)$ and $q_j=y_jx_{j+1}~(j=1,2,\ldots ,n-1)$ segments, and 
$q_n$ an arc in $D^2$ connecting $y_n$ and $x_1$ such that $\bigcup_{i=1}^n(p_i \cup q_j)$ bounds 
the shaded disk in Figure~\ref{fig:disk}. 
Then for an $n$-string link $l$, the knot 
\[l\cup\left(\bigcup_{i=1}^n \Bigl( (p_i\times\{1\})\cup(y_i\times[0,1])\cup (q_i\times\{0\})\Bigr) \right)\]  
is called the {\em closure knot} of $l$. 
Note that the link 
\[l\cup\left(\bigcup_{i=1}^n\Bigl(  p_i\times\{0,1\}\cup y_i\times[0,1] \Bigl) \right)\]
is the {\em closure} $\hat{l}$ of $l$ in the usual sense.

The set of isotopy classes of $n$-string links fixing the endpoints has a monoid structure, 
with composition given by the stacking product and with the trivial $n$-string link ${\bf 1}_n$ as
unit element. Given two $n$-string links $L$ and $L'$, we denote their product by $L\times L'$, 
which is obtained by stacking $L'$ {\em above} $L$ and reparametrizing the ambient cylinder $D^2 \times [0,1]$.

%%%%%%%%%%%%%%%%%%%%
\medskip
\subsection{Clasper}
Clasper is defined by K. Habiro \cite{H}.  
Here we define only tree clasper. 
For a general definition of clasper, we refer the reader to \cite{H}.  

Let $L$ be a (string) link.  
A disk $T$ embedded in $S^3$ (or $D^2\times [0,1]$) is called a {\em tree clasper} for $L$ if it 
satisfies the following three conditions:
\begin{enumerate} 
\item $T$ is decomposed into disks and bands, called {\em edges}, each of which connects two distinct disks.
\item The disks have either 1 or 3 incident edges. We call a disk with 1 incident edge a {\em leaf}.
\item $T$ intersects $L$ transversely, and the intersections are contained in the union of the interiors of the leaves. 
\end{enumerate}
%\end{definition}

Throughout this paper, the drawing convention for claspers are those of \cite[Figure~7]{H}, unless otherwise specified.  

The {\em degree} of a tree clasper $T$ is defined as the number of leaves minus 1.  
A tree clasper of degree $k$ is called a \emph{$C_k$-tree}.  
A tree clasper for a (string) link $L$ is \emph{simple} if each of its leaves intersects $L$ at exactly one point.   
Let $T$ be a simple tree clasper for an $n$-component (string) link $L$.  
The \emph{index} of $T$ is the collection of all integers $i$ such that $T$ intersects the $i$-th component of $L$.

Given a $C_k$-tree $T$ for a (string) link $L$, there is a procedure to construct a
framed link $\gamma(T)$ in a regular neighborhood of $T$. 
{\em Surgery along} T means surgery along $\gamma(T)$. 
Since there exists an orientation-preserving homeomorphism, fixing the boundary, 
from the regular neighborhood $N(T)$ of $T$ to the manifold $N(T)_T$ obtained from $N(T)$ by surgery along $T$, 
surgery along $T$ can be regarded as a local move on $L$. 
We say that the resulting link $L_T$ is obtained from $L$ by surgery along $T$. 
For example, surgery along a simple $C_k$-tree is a local move as illustrated in Figure~\ref{ckmove}.
 
Similarly, for a disjoint union of trees $T_1 \cup \ldots \cup T_m$ for $L$, 
we can define $L_{T_1 \cup \ldots \cup T_m}$ as a link obtained by surgery along $T_1 \cup \ldots \cup T_m$. 
We often regard $L \cup T_1 \cup \ldots \cup T_m$ as $L_{T_1 \cup \ldots \cup T_m}$.

The $C_k$-equivalence is an equivalence relation on (string) links generated by surgeries 
along $C_k$-tree claspers and isotopies.  
We use the notation $L {\sim}_{C_{k}} L'$ for $C_k$-equivalent (string) links $L$ and $L'$.  
%\stackrel

%%%%%%%%%%%%%%%%%%%%%%%%%%%%%%%%
\medskip
\subsection{Linear trees and planarity} 
\begin{figure}[!h]
 \begin{center}
\includegraphics[width=.7\linewidth]{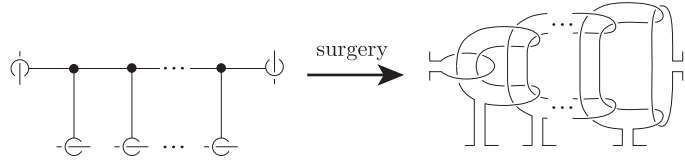}
 \caption{Surgery along a simple tree.} \label{ckmove}
 \end{center}
\end{figure}

For $k\ge 3$, a $C_k$-tree $T$ having the shape of the tree clasper in Figure~\ref{ckmove} is called a \emph{linear} $C_k$-tree. 
The left-most and right-most leaves of $T$ in Figure~\ref{ckmove} are called the \emph{ends} of $T$.

Now suppose that $T$ is a linear $C_k$-tree for some knot $K$, 
and denote its ends by $f$ and $f'$.  
Then the remaining $k-1$ leaves of $T$ can be labeled from $1$ to $k-1$, by 
travelling along the boundary of the disk $T$ from $f$ to $f'$ so that all leaves are visited.  
We say that $T$ is \emph{planar} if, when traveling along $K$ from $f$ to $f'$, either following or 
against the orientation, the labels of the leaves met successively are strictly increasing.  

\medskip
\begin{lemma} (\cite[Lemma 3.2]{MY}) \label{nonplanar}
Let $T$ be a non-planar linear tree clasper for a knot $K$.
Then $P_0(K_T;t)=P_0(K;t)$. 
\end{lemma}

%%%%%%%%%%%%%%%%%%%%%%%%%%%
\medskip
\subsection{Presentation of link-homotopy classes for string links}
\begin{figure}[h]
  \begin{center}
\includegraphics[width=.3\linewidth]{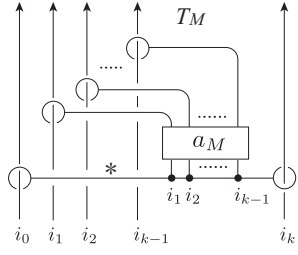}
    \caption{The $C_k$-trees $T_M$ and ${T}^{-1}_M$}
    \label{T_M}
  \end{center}
\end{figure}

Let $\mathcal{M}_k$ denote the set of all sequences $m_0m_1 \ldots m_{k}$ of $k+1$ non-repeating 
integers from $\{ 1,2,\ldots ,n\}$ 
such that $m_0<m_l<m_{k}$ for $1 \leq l \leq k-1$.
Let $i_0 i_1\ldots i_k$ be a subsequence of $12\cdots n$, 
and let $a_M$ be a permutation of $\{i_1,i_2,\ldots ,i_{k-1}\}$.  
Then $M=i_0 a_M(i_1)\ldots a_M(i_{k-1})i_k$ is in $\mathcal{M}_k$ 
and all elements of $\mathcal{M}_k$ can be realized in this way.  
Let $T_M$ be the simple linear $C_k$-tree for ${\bf 1}_n$ as illustrated in Figure~\ref{T_M},
where \fbox{~~~$a_M$~~~} is the unique positive $k$-braid which defined the permutation $a_M$ 
and such that every pair of strings crosses at most one.
In the figure, we also implicitely assume that all edges of $T_M$ \emph{overpass} all 
components of ${\bf 1}_n$. 
Let $T^{-1}_M$ be the $C_k$-tree obtained from $T_M$ by inserting a negative half-twist 
in the $*$ marked edge in Figure~\ref{T_M}.
We remark that a \lq positive' half twist is chosen instead of 
a \lq negative' one  in \cite{MY}. Here we choose negative one for a technical reason 
for the proof of Theorem~\ref{main2}.

Denote respectively by $V_M$ and $V_M^{-1}$ the $n$-string links obtained from ${\bf 1}_n$ 
by surgery along $T_{M}$ and ${T}^{-1}_M$.  
The following theorem is stated  in \cite{MY} as 
a slight modified  version of Theorem~4.3 in \cite{yasuhara}. 

\medskip
\begin{theorem} (\cite[Theorem 4.3]{yasuhara},\cite[Theorem 4.1]{MY}) \label{standerd form}
Any $n$-string link $l$ is link-homotopic to string link $l_1\times l_2 \times \cdots \times l_{n-1}$, where 
\[
  l_i= \prod_{M\in \mathcal{M}_i} V_M^{x_M}\textrm{ , where }\quad
   x_M =\left\{ 
          \begin{array}{ll}
               \mu_l(M) & \text{ if } i=1 \\
               \mu_l(M)-\mu_{l_1 l_2\ldots l_{i-1}}(M) & \text{ if } i \geq 2.  \\
          \end{array} 
            \right.
\]
Here, in the product $\prod_{M\in \mathcal{M}_i} V_M^{x_M}$, 
the string link $V_M^{x_M}$ appears in the lexicographic order of $\mathcal{M}_i$.  
\end{theorem}
%%%%%%%%%%%%%%%%%%%%%%%%%%%%%%%
\medskip
\subsection{The correction term}
Let $K_M^x$ and $K_{M, N}^{x,y}$ denote the knot closures of the string links $V^{x}_M$ and  
 $V_M^x\times V_N^y$ respectively, where  
$x$ and $y$ are integers, and $M$ and $N$ are subsequences of $12\ldots n$.

Let $I=i_1i_2\ldots i_{2k+2}$ be a sequence of $\{1,2,\ldots ,n\}$ without repeated number.
Let $\varphi_I$ be a bijection from $\{1,2,\ldots  ,2k+2\}$ to $\{i_1,i_2,\ldots ,i_{2k+2}\}$ 
which sends any $j$ to $i_j$. 
Let $\mathcal{S}$ be the set of pairs $(M,M')$ such that $M$ and $M'$ are non-successive subsequences 
of $12\ldots (2k+2)$ with length $k+1$, $1<M$ and $\{M\} \cap \{M'\}=\emptyset$.
Then for a link $L$ with vanishing Milnor link-homotopy invariants of length $\leq k$, 
$\delta_L(I)$ is defined by  
\begin{align*}
- \frac{1}{(2k+1)!2^{2k+1}}\sum_{(M,M') \in \mathcal{S}} \left( \log \frac{P_0\left(K_{M,M'}^{\overline{\mu}_{L}
(\varphi_I(M)),\overline{\mu}_{L}(\varphi_I(M'))}\right)}{P_0\left(K_{M}^{\overline{\mu}_{L}
(\varphi_I(M))}\right)P_0\left(K_{M'}^{\overline{\mu}_{L}(\varphi_I(M'))}\right)} \right)^{(2k+1)},
\end{align*} \\
where $\varphi_I(m_1 m_2 \ldots m_{k+1})$ means the sequence $\varphi_I(m_1) \varphi_I(m_2) \ldots 
\varphi_I(m_{k+1})$ for a sequence $m_1  m_2 \ldots m_{k+1}$.  
We note that the Milnor invariants of length $k+1$ for $L$ are integer valued invariants and that 
 they are given by linear combinations of $P_0^{(k)}$'s by \cite[Theorem 1.2]{MY}.  
We also note that $\delta_L(I)$ is a link-homotopy invariant of $L$. 

\medskip
\begin{example}
Let $I=i_1i_2i_3i_4$ be a sequence of $\{1,2,3,4\}$ without repeated number.
Then $\mathcal{S}$ consists of a single pair $(13,24)$ of non-successive subsequences 
of $1234$ and 
\[\delta_L(I)=
- \frac{1}{3!2^{3}} \left( \log \frac{P_0\left(K_{13,24}^{\overline{\mu}_{L}(i_1 i_3),
\overline{\mu}_{L}(i_2 i_4)}\right)}{P_0\left(K_{13}^{\overline{\mu}_{L}(i_1 i_3)}\right)
P_0\left(K_{24}^{\overline{\mu}_{L}(i_2 i_4)}\right)} \right)^{(3)}.
\]
\end{example}

%%%%%%%%%%%%%%%%%%%%%%%%%%
\medskip
\subsection{Calculus of claspers for parallel claspers}\label{sec:parallel}
We shall need  the following lemma  for parallel tree claspers which is given in 
\cite{MY}.  
For a positive integer $m$, an \emph{$m$-parallel} tree  means 
a family of $m$ parallel copies of a tree clasper.

\begin{figure}[h]
  \begin{center}
\includegraphics[width=.9\linewidth]{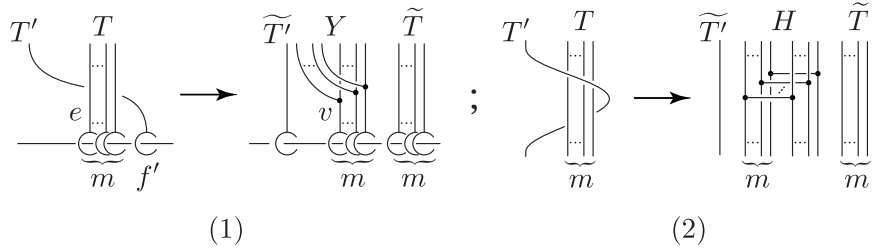}
    \caption{Leaf slide and edge crossing change involving parallel trees}
    \label{slide}
  \end{center}
\end{figure}

\medskip
\begin{lemma}(\cite[Lemma 2.2]{MY}) \label{sliding}
Let $m$ be a positive integer. 
Let $T$ be an $m$-parallel $C_k$-tree for a (string) link $L$,
and $T'$ be a $C_{k'}$-tree for $L$.  
Here $T$ and $T'$ are disjoint.  \\
(1) (Leaf slide)
Let $\widetilde{T}\cup\widetilde{T'} $ be obtained from $T \cup T'$ by sliding a leaf $f'$ of $T'$ 
over $m$ parallel leaves of $T$ (see Figure~\ref{slide}~(1)). 
Then, $L_{T \cup T'}$ is ambient isotopic to $L_{\widetilde{T} \cup \widetilde{T'} \cup Y \cup C }$, 
where $Y$ denotes the $m$ parallel copies of a $C_{k+k'}$-tree obtained by inserting a vertex $v$ in the 
edge $e$ of $T$ and connecting $v$ to the edge incident to $f'$ as shown in Figure~\ref{slide}~(1) and where 
$C$ is a disjoint union of $C_{k+k'+1}$-trees for $L$. \\
(2) (Edge crossing change)
Let $\widetilde{T}\cup\widetilde{T'} $ be obtained from $T \cup T'$ by passing an edge of $T'$ across 
$m$ parallel edges of $T$ (see Figure~\ref{slide}~(2)). 
Then, $L_{T \cup T'}$ is ambient isotopic to $L_{\widetilde{T} \cup \widetilde{T'} \cup H \cup C' }$, 
where $H$ denotes the $m$ parallel copies of a $C_{k+k'+1}$-tree obtained by inserting a vertices in both 
edges, and connecting them by an edge as shown in Figure~\ref{slide}~(2) and where $C'$ is a disjoint union 
of $C_{k+k'+2}$-trees for $L$.      
\end{lemma}

\medskip
\begin{remark}\label{rem:sliding}
Leaf slides between $C_k$-trees for ${\bf 1}_n$ with the same index 
can be realized by link-homotopy, since it 
is realized by surgery along trees intersecting some component 
of ${\bf 1}_n$ more than twice and since a surgery along such trees is realized by link-homotopy \cite[Lemma 1.2]{FY}. 
Hence, in Subsection 2.4,  $V_M^{x_M}$ is link-homotopic to $({\bf 1}_n)_{T_M^{x_M}}$, 
where  $T_M^{x_M}$ is $|x_M|$ parallel copies of $T_M^{x_M/|x_M|}$.  
Note that $V_M^{x_M}=({\bf 1}_n)_{T_M^{x_M}}$ if $x_M\geq 0$.
\end{remark}

%%%%%%%%%%%%%%%%%%%%%%%%%%%%%%%%%%%%%%%%%%%%%%%%%%%%%%%%%%%%%%%%%%%%%%
%%%%%%%%%%%%%%%%%%%%%%%%%%%%%%%%%%%%%%%%%%%%%%%%%%%%%%%%%%%%%%%%%%%%%%
\bigskip
\section{Proof of Theorem \ref{main}}
%%%%%%%%%%%%%%%%%%%%%

Our strategy of the proof is similar to that in \cite[Proof of Theorem 1.1]{MY}. 
In fact, we will use terms \lq good position' and \lq balanced'  for trees 
which are defined in \cite[Proof of Theorem 1.1]{MY} 
and deform, up to $C_n$-equivalence,  a balanced set of trees with keeping it balanced as well. 
The big difference is that we have to treat $C_k$-trees $(n=2k+2)$ while 
they did not need to do. 
We will repeat same arguments as \cite[Proof of Theorem 1.1]{MY} part way. 
We remark that a finite type invariant of degree $\leq n-1$ is an invariant of $C_n$-equivalence \cite{H}, 
in particular $\log(P_0(K))^{(n-1)}$ is an additive invariant of $C_n$-equivalence.

Let $L=\bigcup_{i=1}^{n}L_i$ be an $n$-component link in $S^3$.
Let $I$ be a sequence of $2k+2$ distinct elements of $\{ 1, 2, \ldots , n\}$.
It is sufficient to consider here the case $2k+2=n$, because, if $2k+2<n$, 
we have that $\overline{\mu}_L(I)=\overline{\mu}_{\bigcup_{i \in \{ I \} } L_i}(I)$.
We may further assume that $I=12\ldots n$ without loss of generality.
Indeed, for any permutation $I'$ of $12\ldots n$, we have that $\overline{\mu}_L(I')=\overline{\mu}_{L'}(12\ldots n)$, 
where $L'$ is obtained from $L$ by reordering the components appropriately.

\begin{figure}[h]
\begin{center}
\includegraphics[width=.2\linewidth]{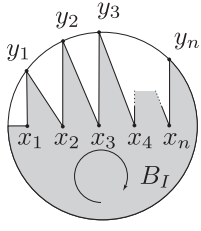}
\end{center}
\caption{The $2n$-gon $B_I$ lying in the unit disk $D^2$. }\label{fig:disk}
\end{figure}

Let $B_I$ be an $I$-fusion disk for $L$.  
Up to isotopy, we may assume that the $2n$-gon $B_I$ lies in the unit disk $D^2$ as shown in Figure~\ref{fig:disk}, 
where the edges $p_j$ ($j=1,2, \ldots ,n$) are defined by $p_j=x_jy_j$.  
We may furthermore assume that $L\cup B_I$ lies in the cylinder $D^2\times [0,1]$,   
such that $B_I\subset (D^2\times \{ 0 \})$, and such that
$$L\cap \partial (D^2\times [0,1]) 
= \bigcup_{j=1}^n \Bigl( (p_j\times \{0,1\})\cup (\{ y_j \}\times [0,1]) \Bigr). $$ 

\noindent 
Then, we obtain an $n$-string link $\sigma$ whose closure $\hat{\sigma}$ is the link $L$, by setting 
%\begin{equation}\label{eq:defsigma}
 \[\sigma := \overline{L\setminus (L\cap \partial (D^2\times [0,1]))}. \]
%\end{equation}

For an $n$-string link $\sigma=\bigcup_{i=1}^n \sigma_i$ and for a subsequence $J$ of $I=12\ldots n$, 
we denote by $\sigma(J)$ the knot 
\[\overline{ \biggl( \Bigl( \bigcup_{j\in \{J\}} \hat{\sigma}_j \Bigr) \cup \partial B_I \biggr) 
\setminus \biggl( \Bigl( \bigcup_{j\in \{J\}} \hat{\sigma}_j \Bigr) \cap \partial B_I \biggr) 
}. \]
We note that the knot $\sigma(I)$ is equal to the closure knot of $\sigma$ defined in Subsection~2.1.
%%%%%%%%%%%%%%%%%%%%%

Let $\sigma$ be the $n$-string link with the closure $L$ defined as above. 
By combining Theorem~\ref{standerd form}, Remark~\ref{rem:sliding}, and the assumption that Milnor 
link-homotopy invariants of length $\leq k$ vanish,
$\sigma$ is link-homotopic to $l_k \times \cdots \times l_{2k+1}$, where $l_i=\prod_{M \in \mathcal{M}_i}
 ({\bf 1}_n)_{T_M^{x_M}}$.
Therefore there is a disjoint union $R_1$ of simple $C_1$-trees whose leaves intersect a single component of 
$l_k \times \cdots \times l_{2k+1}$ such that 
$R_1$ is disjoint from $\bigcup_{i=k}^{2k+1}\bigcup_{M \in \mathcal{M}_i} T_M^{x_M}$ and
\begin{align*}
    \sigma = (l_k \times \cdots \times l_{2k+1})_{R_1}.
\end{align*}

\noindent
Set   
\[G := \bigcup _{i=k}^{2k+1}\left(\bigcup _{M \in \mathcal{M}_i} T_M^{x_M}\right), \]
then we have 
\[ \sigma = ({\bf 1}_n)_{G \cup R_1}. \]

\noindent
Moreover, by combining Lemmas 4.3 and 4.4 in \cite{MY}, we have the following: 
\begin{align}\label{claim:milnor1}
\overline{\mu}_L(I) \equiv x_I ~\pmod{\Delta _L(I)}, %~\text{and}
\end{align}
and
\begin{align}\label{claim:milnor2}
\Delta _L(I)= \gcd \{ x_M \mid M < I, M\neq I\}.
\end{align}

A tree for ${\bf 1}_n$ is said to be in {\em good position} if each component of ${\bf 1}_n$ 
underpasses all edges of the tree.
Note that each tree of $G$ is in good position.
On the other hand, a tree of $R_1$ may not be in good position.  
We now replace $R_1$ with some trees with good position up to $C_n$-equivalence. 
By \cite[Proposition 4.5]{H}, we have
\[
    ({\bf 1}_n)_{G \cup R_1} \sim _{C_n} ({\bf 1}_n)_{G \cup R},
\]
where $R$ is a disjoint union of simple trees for ${\bf 1}_n$ in good position and intersecting some component 
of ${\bf 1}_n$ more than once.

It follows from Lemma \ref{nonplanar} that for any $J < I$,  
\[
    P_0(({\bf 1}_n)_{G \cup R}(J)) = P_0(({\bf 1}_n)_{\widetilde{G} \cup R}(J)),
\]
where $\widetilde{G}$ is obtained from $G$ by eliminating non-planar trees for $({\bf 1}_n)(I)$.
That is, 
\[
    \widetilde{G} = \bigcup _{i=k}^{2k+1}\left(\bigcup _{M \in \mathcal{M}_i, M<I} T_M^{x_M}\right).
\]  
Here, since $\Delta_L(I)$ divides all $x_M$ with $M<I$ and $M\neq I$ by  (\ref{claim:milnor2}), we assume that 
each $T_M^{x_M}$ is a disjoint union of $\Delta_L(I)$ parallel copies of $T_M^{x_M/|x_M|}$. 

We now define the {\em weight} of a tree $t$ for the trivial knot as a subset of $\{ 1, 2, \ldots ,n\}$ 
and denote it by $w(t)$. 
A disjoint union $g_1 \cup \ldots \cup g_s$ of $s$ trees (possibly parallel) for the trivial knot 
$U$ is {\em balanced} if each tree has a weight such that 
\begin{align*}
    ({\bf 1}_n)_{\widetilde{G} \cup R}(J) \sim_{C_n} U_{\left(\bigcup_{w(g_i)\subset \{ J\}}g_i\right)},
\end{align*}
for any $J<I$.

For the knot $({\bf 1}_n)_{\widetilde{G} \cup R}(I)$, we may think of  
$\widetilde{G}\cup R$ as trees $F$ for the trivial knot $U(=({\bf 1}_n)(I))$.  
We assign the index of each tree of $F$ as weight. 
Here we recall that the index of a tree for a (string) link is the collection of all integers $i$ such 
that the tree intersects the $i$-th component of the (string) link. 
We may assume that each tree with index $\subset\{J\}$ is also a tree for $({\bf 1}_n)(J)$. 
Then it is obvious that  for any $J<I$
\[({\bf 1}_n)_{\widetilde{G} \cup R}(J) =(({\bf 1}_n)(J))_{\left(\bigcup_{w(g)\subset \{ J\}}g\right)} ,\]
where $\bigcup_{w(g)\subset \{ J\}}g$ means the union of trees $g$ of $F$ with weight $\subset \{J\}$. 
Since $\widetilde{G} \cup R$ is in good position, 
$(({\bf 1}_n)(J))_{\left(\bigcup_{w(g)\subset \{ J\}}g\right)}$ and 
$U_{\left(\bigcup_{w(g)\subset \{ J\}}g\right)}$  
have a common diagram in $D^2\times\{0\}$, 
and hence they are ambient isotopic. 
In particular $F$ is balanced. 

%Since $\widetilde{G} \cup R$ is in good position, $F$ is balanced if we assign the index of 
%each tree of $F$ as weight. 

\medskip
\begin{remark}\label{rem:weight}
When we perform a leaf slide or an edge crossing change between two trees 
in a balanced union of trees as in Lemma \ref{sliding}, 
we assign the union of weights as weight to each of new trees. 
More precisely, in Lemma \ref{sliding} (1) (resp.~(2)), we assign the weights $w(T)$ and $w(T')$ to $\widetilde{T}$ 
and $\widetilde{T'}$ respectively, and assign the union $w(T)\cup w(T')$ to $Y$ (resp.~$H$) 
and each connected component of $C$ (resp.~$C'$). 
We note that the union of resulting trees is also balanced. 
\end{remark}

So far, the proof is the same as  \cite[Proof of Theorem 1.1]{MY}. 
In  \cite{MY}, they deform $F$ into a balanced union of \lq localized' tree for the trivial knot $U$ up to $C_n$-equivalence. 
But in that case, there are no $C_k$-trees $(n=2k+2)$. 
The main difficulty of our proof is how to treat such $C_k$-trees. 
In the following, we first deform $F$ into a balanced union of \lq separated' trees for $U$ except 
for $C_k$-trees, and then deform the $C_k$-trees into suitable shape. 
Here \lq localized' implies \lq separated'. 
We use \lq separated' instead of \lq localized', since we notice that we do not need to 
such a strong condition as \lq localized'. So we also slightly modify \cite[Proof of Theorem 1.1]{MY} 
in this sense.  

For ${M\in{\mathcal{M}}_i}~(i=k,\ldots ,2k+1)$, we denote by $t_M^{\pm1}$ the tree of $F$ which corresponds to 
$T_M^{\pm1}$ of $\widetilde{G}$. Set
\[d:=\bigcup_{M\in{\mathcal M}_k,M<I}t^{x_M}_M,~\text{and}~D:=U_d.\]
Then $U_F$ is obtained from $D$ by surgery along the trees of $F\setminus d$. 
By using leaf slides and edge crossing changes, we will deform, up to $C_n$-equivalence, $F$ 
into a balanced set of \lq separated' trees for $U$ with fixing $d$ 
as in claim below.

\medskip
\begin{claim}\label{claim1}
The knot $U_F$ is $C_n$-equivalent to the connected sum 
$D\# (\#_{J< I} C^J)$ of knots $D$ and $C^J~(J< I)$, where 
$C^J$ is a knot obtained from the trivial knot by surgery along 
a disjoint union $F_J$ of trees with weight $\{J\}$ and the set 
$d\cup(\bigcup_{J<I} F_J)$ is balanced. 
Moreover $F_I$ consists of the parallel tree $t_{I}^{x_I}$ and 
$\Delta_L(I)$-parallel $C_{2k+1}$-trees. 
\end{claim}

We define that a tree has {\em full weight} if the degree of the tree plus 
one is equal to the number of weight of the tree. 
We define that a tree is {\em repeated} if the degree of the tree 
is more than or equal to the number of weight of the tree. 

\medskip
\begin{proof}
We take mutually disjoint 3-balls $N_J~(\{J\}\in2^{\{1,2,\ldots ,n\}})$ such that  
$(N_J,N_J\cap D)$ is a trivial ball-arc pair and $N_J\cap d=\emptyset$.
By using leaf slides and edge crossing changes, i.e., by Lemma~\ref{sliding} and Remark~\ref{rem:weight}, 
we may assume that all trees except for $d$ with weight $\{J\}$ are contained in the interior of $N_J$ 
up to $C_n$-equivalence with keeping the set of trees balanced.  
Then we have that $U_F$ is $C_n$-equivalent to 
$D\#(\#_{J<I}C^J)$ and $C^J$ is obtained from $U$ by 
surgery along trees contained in $N_J$, which are trees with weight $\{J\}$. 

To complete the proof, we need to show that the trees $F_I$ in $N_I$ consists of 
the parallel tree $t_I^{x_I}$ and some $\Delta_L(I)$-parallel  $C_{2k+1}$-trees. 
Since $t_I$ is $C_{2k+1}$-tree and $n=2k+2$, by Lemma~\ref{sliding} we 
can freely move $t_I^{x_I}$ into $N_I$ up to 
$C_{n}$-equivalence. 
By Remark~\ref{rem:weight} and the observation below, 
we see that whenever we apply Lemma~\ref{sliding},  
the new trees we get are repeated or have full weight.  
Moreover trees have full weight only if they are $\Delta_L(I)$-parallel trees. 
Hence we obtain the claim.
\end{proof}
%%%%%%

\medskip
\begin{observation}
We always move $\Delta_L(I)$-parallel trees together. 
If a leaf of new tree obtained by a leaf slide or an edge crossing change interrupts a parallel leaf 
of a parallel tree, then we sweep the new leaf out of the parallel leaf up to $C_n$-equivalence.      
Since the degrees of parallel trees are at least $k$ and the new tree at least $k+1$,  
we can do such sweeping out easily up to $C_{n}$-equivalence by Lemma~\ref{sliding}.

We consider a leaf slide between a full weight $\Delta_L(I)$-parallel tree $t$ and a repeated tree $t'$.
Let $m$ be the degree of $t$ and $l$ the degree of $t'$.   
If $w(t) \cap w(t') = \emptyset$, then a new $C_{m+l}$-tree, 
which is a $\Delta_L(I)$-parallel tree, has a weight 
consisting of at most $m+l+1$ elements and new $C_{m+l+1}$-trees are repeated.
If $w(t) \cap w(t') \neq \emptyset$, then all new trees are repeated.  
      
We consider a leaf slide between full weight parallel trees $t$ and $t'$.
We may assume that the degree of $t$ is at least $k+1$ and the degree of $t'$ is 
at least $k$.  
Then the new trees  are $\Delta_L(I)$-parallel trees with degree at least $n-1$.

A leaf slide between repeated trees and  an edge crossing change for any case   
give only repeated trees. 
\end{observation}

Now we consider $D$ in Claim~\ref{claim1}. 
Let $\mathcal{S}_k^0$ be the set of pairs $(M,M')$ such that $M$ and $M'$ 
are  subsequences of $I$ with length $k+1$, $1<M$, and $\{M\} \cap \{M'\}=\emptyset$.
We also denote by $\mathcal{S}_k$ the subset of $\mathcal{S}_k^0$ such that both sequences $M$ and $M'$ 
are not successive.
We note that 
\begin{align*}
    \bigcup_{(M,M') \in \mathcal{S}_k^0} (t_M^{x_M} \cup t_{M'}^{x_{M'}}) = \bigcup _{M \in \mathcal{M}_k, M<I} t_M^{x_M}=d.
\end{align*}

We separate $d$ into pairwise trees $t_M^{x_M} \cup t_{M'}^{x_{M'}}~((M,M')\in\mathcal{S}_k^0)$   
by leaf slides and edge crossing changes  
between different pair of parallel trees. 
For two parallel trees $t_M^{x_M}$ and $t_N^{x_N}$ which are not pair, 
we note that $w(t_M) \cap w(t_{N})=\{M\}\cap \{N\} \neq \emptyset$. 
Therefore, when we apply leaf slides or edge crossing changes between 
$t_M^{x_M}$ and $t_N^{x_N}$, we obtain new trees with  
degree at least $2k$ which are $\Delta_L(I)$-parallel trees with 
full weight and/or repeated trees.  

%For a tree clasper $t$ for $U$ and a positive integer $x$, we denote by $\overline{t^x}$ a disjoint union of trees 
%obtained from a $x$-parallel tree $t^{x}$ by inserting a negative half twist in an edge of each component 
%of $t^{x}$ so that each component of $\overline{t^{x}}$ is equal to $\overline{t}$.
We denote by $d_M$ 
the parallel tree $t_M^{x_M}$  if $x_M\geq 0$ and 
the disjoint union of trees obtained from $t_M^{|x_M|}$ by 
 inserting a negative half twist in an edge of each component 
so that each component of $d_M$ is equal to $t_M^{-1}$ 
if $x_M< 0$.
Then we note that
\begin{align*}
K_{M}^{x_M}=U_{d_M } \text{ and } K_{M,M'}^{x_M,x_{M'}} = U_{d_M\cup d_{M'}}.
\end{align*}

\noindent
By using  leaf slides, we have 
\begin{align*}
U_{t_M^{x_M} } \sim_{C_{2k}} U_{d_M}  \text{ and } U_{t_M^{x_M} \cup t_{M'}^{x_{M'}}} \sim_{C_{2k}} U_{d_M\cup d_{M'}},
\end{align*}
where the $C_{2k}$-equivalence is realized by surgery along repeated $C_{2k}$-trees. 
Hence we have that $D$ is $C_n$-equivalent to the connected sum 
$D'\#(\#_{(M,M')\in\mathcal{S}_k^0}U_{d_M\cup d_{M'}})$, 
where $D'$ is obtained from $U$ by surgery along a union of repeated trees. 
Set 
\[\widetilde{D}:=\#_{(M,M')\in\mathcal{S}_k^0}U_{d_M\cup d_{M'}},~\text{and}~  
\widetilde{d}:=\bigcup _{M \in \mathcal{M}_k, M<I} {d_M}.\]

Hence $D\# (\#_{J < I} C^J)$ is $C_n$-equivalent to 
$D'\#\widetilde{D}\# (\#_{J < I} C^J)$. 
By the same reason as Claim~\ref{claim1}, we have the following claim. 

\medskip
\begin{claim}\label{claim2}
 The knot $U_F$ is $C_n$-equivalent to 
$\widetilde{D}\# (\#_{J < I} \widetilde{C}^J)$, where 
$\widetilde{C}^J$ is a knot obtained from the trivial knot by surgery along 
a disjoint union $\widetilde{F}_J$ of trees with weight $\{J\}$ and 
$\widetilde{d}\cup(\bigcup_{J<I}\widetilde{F}_J)$ is balanced. 
Moreover $\widetilde{F}_I$ consists of the parallel tree $t_{I}^{x_I}$ and 
$\Delta_L(I)$-parallel  $C_{2k+1}$-trees. 
\end{claim}

Now we have 
\[\begin{array}{rcl}
\displaystyle\sum_{J<I}(-1)^{|J|}(\log P_0(L_J))^{(n-1)}
&=&\displaystyle\sum_{J<I}(-1)^{|J|}
(\log P_0(\widetilde{D}^J\#(\#_{J'<J}\widetilde{C}^{J'})))^{(n-1)}\\
&=&\displaystyle\sum_{J<I}(-1)^{|J|}
(\log P_0(\widetilde{D}^J))^{(n-1)}\\
&&\displaystyle+\sum_{J<I}(-1)^{|J|}\sum_{J'<J}(\log P_0(\widetilde{C}^{J'}))^{(n-1)},
\end{array}
\]
where $\widetilde{D}^J$ is a knot obtained from $U$ by surgery along 
the union of trees in $\widetilde{d}$ whose weights are subsets of $\{J\}$.  
Note that $\widetilde{D^I}=\widetilde{D}(=\#_{(M,M')\in\mathcal{S}_k^0}K_{M,M'}^{x_M,x_{M'}})$. 

For $J'<I$ and $J'\neq I$, the coefficient of $(\log P_0(\widetilde{C}^{J'}))^{(n-1)}$ in 
\[\displaystyle\sum_{J<I}(-1)^{|J|}\sum_{J'<J}(\log P_0(\widetilde{C}^{J'}))^{(n-1)}\]
is equal to 0, since 
\[\sum_{J' < J <I}(-1)^{|J|}=
\sum_{\{J'\}\subset\{J\}\subset\{I\}\setminus\{a\}}(-1)^{|\{J\}|}+
\sum_{\{J'\}\subset\{J\}\subset\{I\}\setminus\{a\}}(-1)^{|\{J\}\cup\{a\}|}=0,
\]
where $a$ is an element in $\{I\}\setminus\{J'\}$. 
Hence we have 
\begin{align}\label{eq1}
\begin{array}{l}
\displaystyle\sum_{J<I}(-1)^{|J|}(\log P_0(L_J))^{(n-1)}\\
\hspace*{3em}=\displaystyle\sum_{J<I}(-1)^{|J|}(\log P_0(\widetilde{D}^J))^{(n-1)}
+(\log P_0(\widetilde{C}^{I}))^{(n-1)}. 
\end{array}
\end{align}
Note that $(-1)^n=(-1)^{2k+2}=1$.

On the other hand, 
\begin{align*}
   &\sum_{J<I} (-1)^{|J|}(\log P_0(\widetilde{D}^J))^{(n-1)} \\
   &\hspace*{3em}=\sum_{J<I;J\neq I} (-1)^{|J|}(\log P_0(\widetilde{D}^J))^{(n-1)} +
(\log P_0(\widetilde{D}^I))^{(n-1)}\\
&\hspace*{3em} =\sum_{J<I;J\neq I} (-1)^{|J|}(\log P_0(\widetilde{D}^J))^{(n-1)} +
\sum_{(M,M')\in {\mathcal S}_k^0}(\log P_0(K_{M,M'}^{x_M,x_{M'}}))^{(n-1)}\\
\end{align*}
For a subsequence $M$ of $I$ with length $k+1$, the coefficient of 
$(\log P_0(K_M^{x_M}))^{(n-1)}$ in $\sum_{J<I;J\neq I} (-1)^{|J|}(\log P_0(\widetilde{D}^J))^{(n-1)}$ is 
\[
\sum_{M<J; J\neq I} (-1)^{|J|} =\sum_{i=0}^k 
\left(\begin{array}{c}
{k+1} \\
i
\end{array}\right)\times
 (-1)^{i+k+1} = -1.
\]
This implies that 
\begin{align*}
   & \sum_{J<I} (-1)^{|J|}(\log P_0(\widetilde{D}^J))^{(n-1)} \\
&  =\sum_{(M,M')\in {\mathcal S}_k^0} \Bigl( (\log P_0(K_{M,M'}^{x_M,x_{M'}}))^{(n-1)} - (\log P_0(K_{M}^{x_M}))^{(n-1)} - 
(\log P_0(K_{M'}^{x_{M'}}))^{(n-1)} \Bigr)\\
&=\sum_{(M,M')\in {\mathcal S}_k^0}\Biggl( \log \frac{P_0(K_{M,M'}^{x_M,x_{M'}})}{P_0(K_{M}^{x_M})P_0(K_{M'}^{x_{M'}})} \Biggr)^{(n-1)}.
\end{align*}
If $(M,M') \in \mathcal{S}_k^0 \setminus  \mathcal{S}_k$, then $\widetilde{d_M}$ and 
$\widetilde{d_{M'}}$ are separated by a 2-sphere since either  $M$ or $M'$ is 
a successive sequence. Hence we have 
\[
  K_{M,M'}^{x_M,x_{M'}} = K_{M}^{x_M} \#  K_{M'}^{x_{M'}}~~((M,M') \in \mathcal{S}_k^0 \setminus  \mathcal{S}_k). 
\]
It follows that 
\begin{align}\label{eq2}
\begin{array}{rcl}
\displaystyle\sum_{J<I} (-1)^{|J|}(\log P_0(\widetilde{D}^J))^{(n-1)}
&=&
\displaystyle\sum_{(M,M') \in \mathcal{S}_k} \Biggl( \log \frac{P_0(K_{M,M'}^{x_M,x_{M'}})}
{P_0(K_{M}^{x_M})P_0(K_{M'}^{x_{M'}})} \Biggr)^{(n-1)}\\
&&\\
&=&-(n-1)!2^{n-1}\delta _L (I). 
\end{array}
\end{align}

We now consider $\widetilde{C}^I$. 
Let $h_1^{\Delta_L(I)}, h_2^{\Delta_L(I)},\ldots,h_r^{\Delta_L(I)}$ be the 
$\Delta_L(I)$-parallel $C_{2k+1}$-trees in $\widetilde{F}_I\setminus t_I^{x_I}$.
Then by using leaf slides and edge crossing changes, we have that 
\[\widetilde{C}^I\sim_{C_n} (|x_I|\times U_{t_I^{x_I/|x_I|}})
\#(\Delta_L(I)\times(\#_{i=1}^r U_{h_i})),\]
where for positive integer $x$ and for a knot $K$, $x\times K$ denotes the connected 
sum of $x$ copies of $K$. 
By combining \cite[Lemma 3.1 and Claim 5.3 (2)]{MY} and (\ref{claim:milnor1}), we have 
\begin{align}\label{eq3}
    - \frac{1}{(n-1)!2^{n-1}} (\log P_0(\widetilde{C}^I))^{(n-1)} 
   \equiv x_I %%%\pmod {\Delta _L(I)}  \notag \\
   \equiv  \overline{\mu}_{L}(I)  \pmod {\Delta _L(I)}.  
\end{align}

It follows from Equations\ (\ref{eq1}), (\ref{eq2}) and (\ref{eq3}) that 
\begin{align*}
 \overline{\mu}_{L}(I) \equiv & - \frac{1}{(n-1)!2^{n-1}} \sum_{J<I} (-1)^{|J|}(\log P_0(L_J))^{(n-1)}
 - \delta _L (I) \pmod {\Delta _L(I)}.
\end{align*} \\

%%%%%%%%%%%%%%%%%%%%%%%%%%%%%%%%%%%%%%%%%%%%%%%%%%%%%%%%%%%%%%%%%%%%%%
%%%%%%%%%%%%%%%%%%%%%%%%%%%%%%%%%%%%%%%%%%%%%%%%%%%%%%%%%%%%%%%%%%%%%%
\bigskip
\section{Proof of Theorem \ref{main2}}

\medskip
\subsection{HOMFLYPT polynomial}

First of all, we recall the definition of the HOMFLYPT polynomial, and mention a few useful
properties. 

The \emph{HOMFLYPT polynomial} $P(L;t,z)\in {\Bbb Z}[t^{\pm 1},z^{\pm 1}]$ of an oriented link $L$ 
is defined by the following formulas:
\begin{enumerate}
\item $P(U;t,z) = 1$, and 
\item $t^{-1}P(L_+;t,z) - tP(L_- ;t,z) = zP(L_0 ;t,z)$, 
\end{enumerate}
where $U$ denotes the trivial knot and where $L_+$, $L_-$ and $L_0$ are three links that are identical 
except in a $3$-ball, where they look as follows: 
\[L_+=\begin{array}{c}
\includegraphics[width=.07\linewidth]{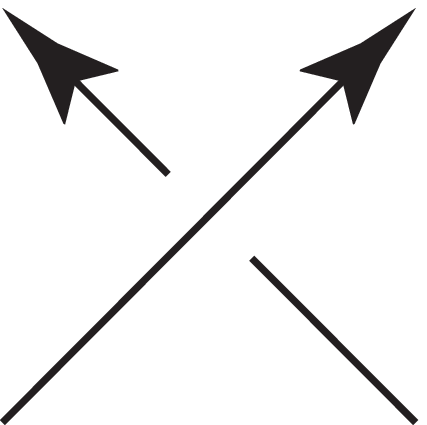}\end{array}~~;~~~~~~
L_+=\begin{array}{c}
\includegraphics[width=.07\linewidth]{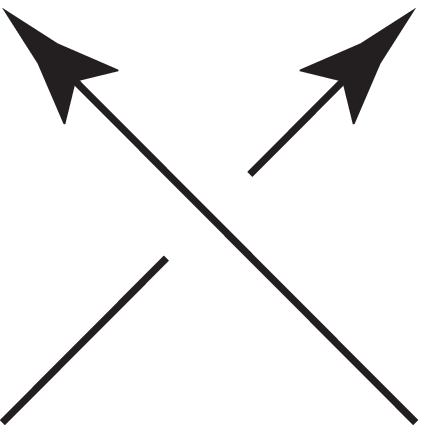}\end{array}~~;~~~~~~
L_+=\begin{array}{c}
\includegraphics[width=.07\linewidth]{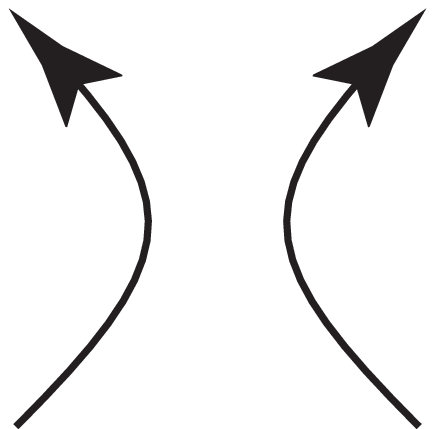}\end{array}~~.
\]
In particular, the HOMFLYPT polynomial of an $r$-component link $L$ is of the form  
 $$ P(L;t,z)=\sum_{k=1}^N P_{2k-1-r}(L;t)z^{2k-1-r}, $$
where $P_{2k-1-r}(L;t)\in {\Bbb Z}[t^{\pm 1}]$ is called the $(2k-1-r)$-th coefficient polynomial of $L$.  
Furthermore, the lowest degree coefficient polynomial of $L$ is given by 
\begin{equation}\label{eq:LM}
 P_{1-r}(L;t)=t^{2Lk(L)}(t^{-1}-t)^{r-1}\prod_{i=1}^r P_0(L_i;t), 
\end{equation}
where $L_i$ is the $i$-th component of $L$ and $Lk(L)$ 
is the total linking numbers, 
see \cite[Prop. 22]{LM}.

%%%%%%%%%%%%%%%%%%%%%%%%%%%%%%%%%%%%%%%%
\medskip
\subsection{Proof  of Theorem \ref{main2}}

\begin{figure}[h]
  \begin{center}
     \includegraphics[scale=.4]{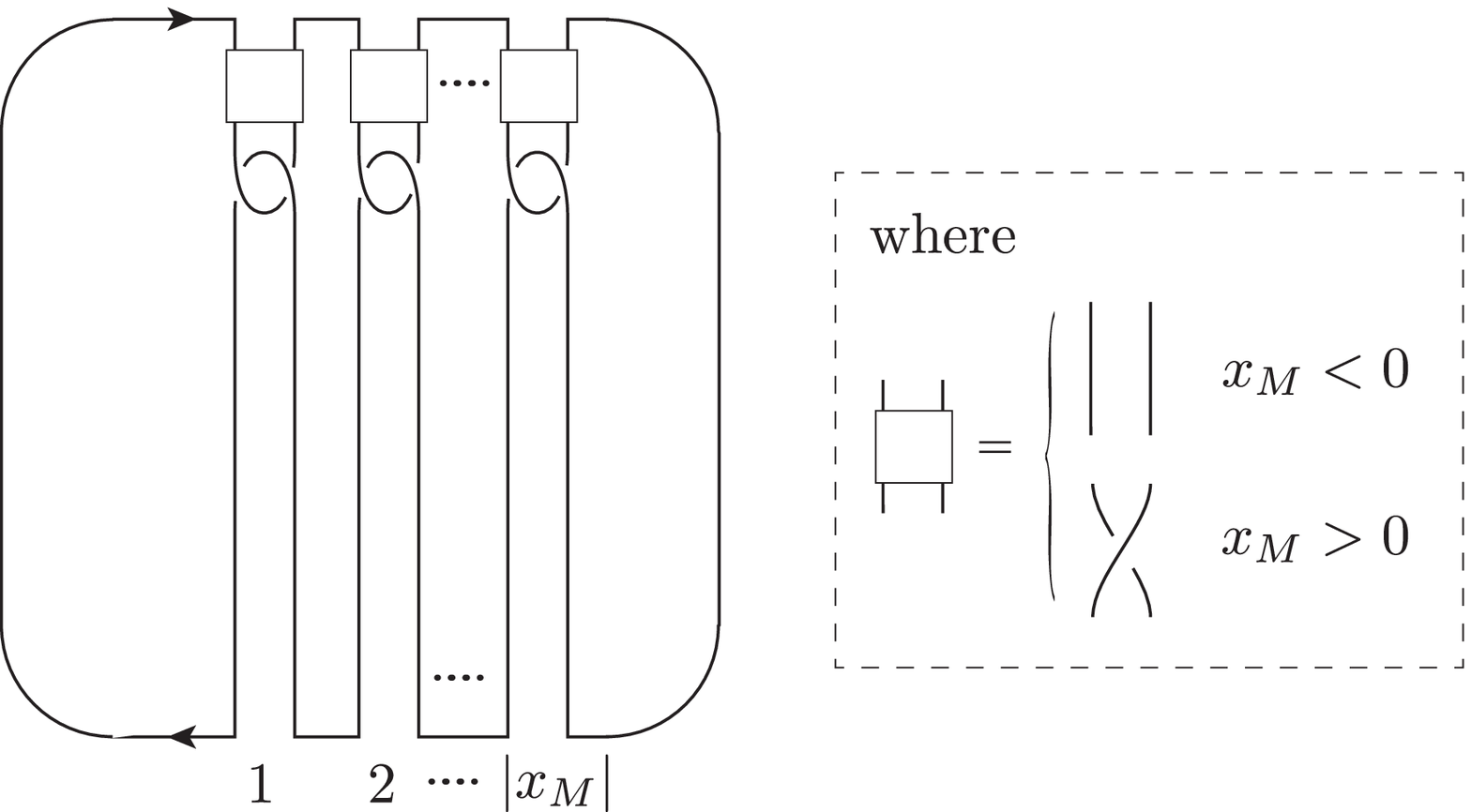}
    \caption{}
    \label{KxM}
  \end{center}
\end{figure}

We may assume that $I=1234$ by the same reason as those in the proof of Theorem~\ref{main}. 
By leaf slides and edge crossing changes, we deform the shape of a disjoint union $d_M$ of $C_1$-trees  
which appears in the proof of Theorem~\ref{main} (here $k=1$)  
so that the knot $U_{d_M}$ 
is as illustrated in Figure~\ref{KxM}, which is ambient isotopic to the trivial knot.  
Since these deformation can be realized by surgery along repeated trees, 
we obtain Theorem~\ref{main} for the case when $k=1$ but different correction term. 
We remark that the difference of correction terms vanishes modulo $\Delta_L(I)$. 
Here we have that the new correction term is 
\begin{align*}
(\log P_0(K(x_{13},x_{24})))^{(3)},
\end{align*}
where $x_{ij}=\overline{\mu}_L(ij)$ and $K(m,n)$ is a knot as illustrated in Figure~\ref{K(m,n)}.
Since $(\log P_0)^{(3)} = P_0^{(3)}$, we have 
\begin{align}\label{new correction term}
  \overline{\mu}_{L}(I) \equiv & - \frac{1}{48} \sum_{J<I} (-1)^{|J|} P_0^{(3)}(L_J) 
+\frac{1}{48} P_0^{(3)}K(x_{13},x_{24}) \pmod {\Delta _L(I)}.
\end{align}

\begin{figure}[h]
  \begin{center}
     \includegraphics[scale=.4]{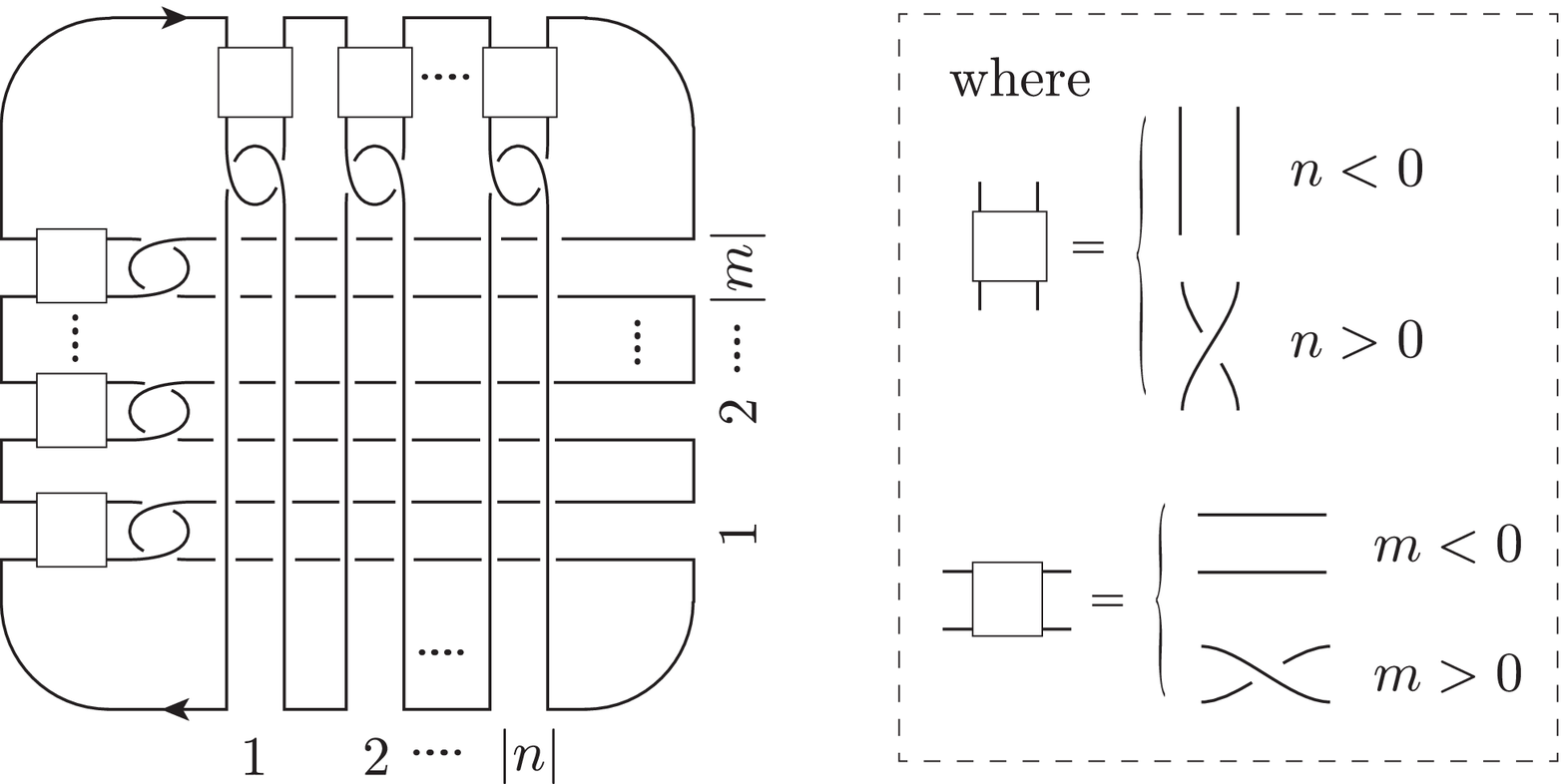}
    \caption{$K(m,n)$}
    \label{K(m,n)}
  \end{center}
\end{figure}

\begin{figure}[h]
  \begin{center}
     \includegraphics[scale=.4]{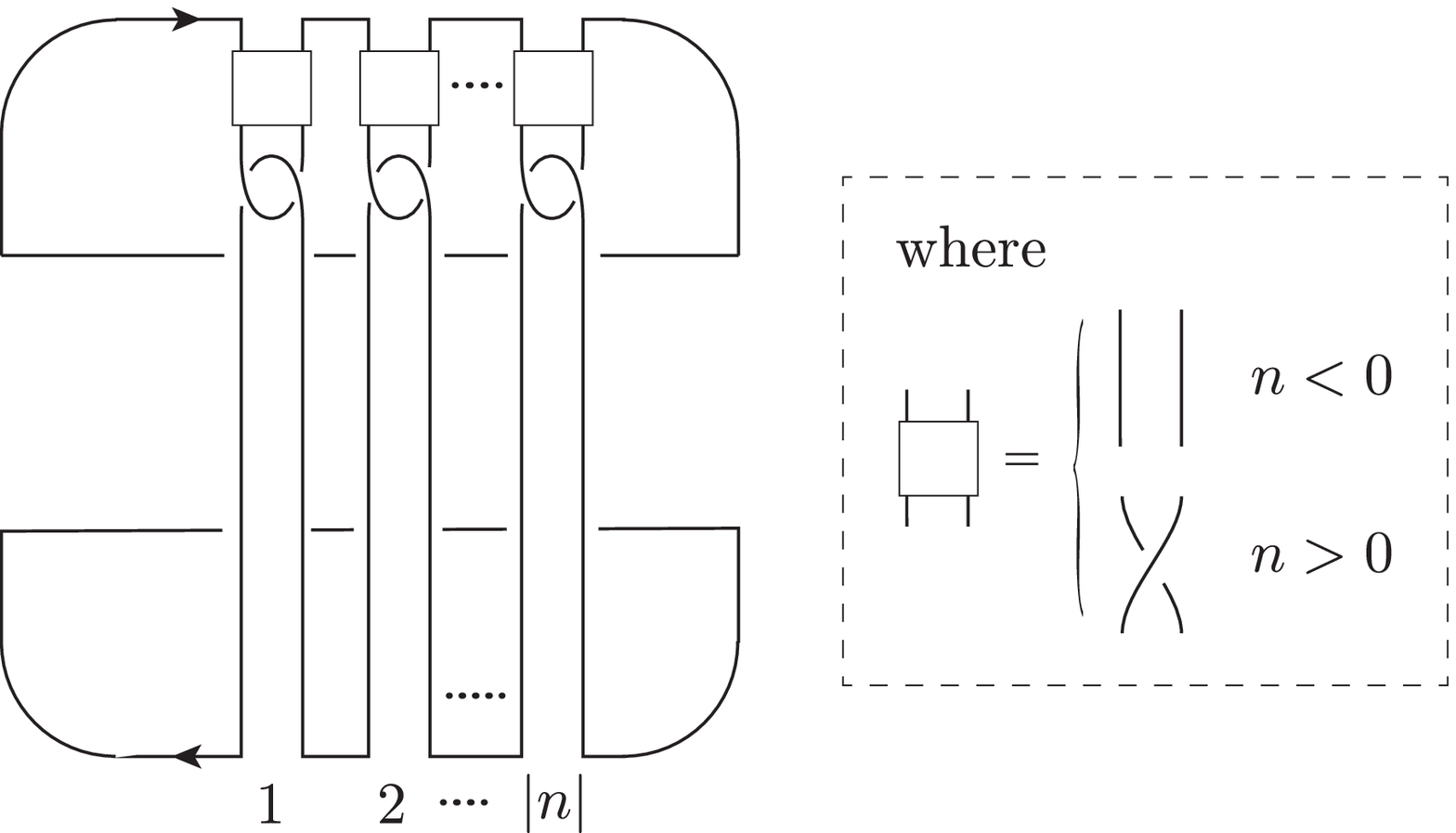}
    \caption{$L(n)$}
    \label{L(n)}
  \end{center}
\end{figure}

We calculate $P_0(K(m,n))$.
Using the relation of the HOMFLYPT polynomial, we obtain the relation
\begin{align}\label{relation of K(m,n)}
   P_0(K(m,n)) = t^{2\varepsilon} P_0(K(m-\varepsilon,n)) + \varepsilon t^{\varepsilon} P_{-1}(L(n)), 
\end{align}
where $L(n)$ is illustrated in Figure~\ref{L(n)}, and $\varepsilon=1$ (resp. $-1$)  
if $m>0$ (resp. $m<0$).   
Since $Lk(L(n))=n$ and each component of $L(n)$ is trivial, it follows from (\ref{eq:LM}) that  
\begin{equation}\label{eq:LM2}
 P_{-1}(L(n);t)=t^{2n}(t^{-1}-t). 
\end{equation}
By combining (\ref{relation of K(m,n)}) and (\ref{eq:LM2}), 
\begin{align*}
   P_0(K(m,n)) = t^{2\varepsilon} P_0(K(m-\varepsilon,n)) + \varepsilon 
t^{2n-1+\varepsilon}-\varepsilon t^{2n+1+\varepsilon}. 
\end{align*}
Since for each $\varepsilon(\in\{-1,1\}$) 
\[\varepsilon 
t^{2n-1+\varepsilon}-\varepsilon t^{2n+1+\varepsilon}=
t^{2n}-t^{2n+2\varepsilon},\]
we have 
\[   P_0(K(m,n))-t^{2n} = t^{2\varepsilon}( P_0(K(m-\varepsilon,n) -t^{2n}),\] 
and hence
\[   P_0(K(m,n))-t^{2n} = t^{2\varepsilon |m|}( P_0(K(0,n) -t^{2n})=t^{2m}( 1 -t^{2n}).\] 
It follows that we have
\begin{align*}
   P_0(K(m,n))= t^{2m} + t^{2n} - t^{2m+2n}, 
\end{align*}
and so we have
\begin{align}\label{P(K(m,n))}
   P_0^{(3)}(K(m,n))
%= 2m(2m-1)(2m-2) + 2n(2n-1)(2n-2) - (2m+2n)(2m+2n-1)(2m+2n-2)  \\
= -24mn(m+n-1). 
\end{align}

By  (\ref{new correction term}) and (\ref{P(K(m,n))}) we have the required formula
\begin{align*}
\overline{\mu}_{L}(I) \equiv - \frac{1}{3!2^{3}}  \sum_{J<I} (-1)^{|J|}P_0^{(3)}(L_J)
-\frac{1}{2} x_{13}x_{24}(x_{13} +x_{24} -1)  \pmod {\Delta _L(I)}.  
\end{align*}

%%%%%%%%%%%%%%%%%%%%%%%%%%%%%%%%%%%%%%%%%%%%%%%%%%%%%%%%%%%%%%%%%%%%%%
%%%%%%%%%%%%%%%%%%%%%%%%%%%%%%%%%%%%%%%%%%%%%%%%%%%%%%%%%%%%%%%%%%%%%%

\end{document}